\documentclass{equadiff}

\usepackage{amsmath,amssymb,enumitem}
\allowdisplaybreaks
\let\pa\partial  
\let\na\nabla  
\let\eps\varepsilon  
\newcommand{\N}{{\mathbb N}}  
\newcommand{\R}{{\mathbb R}} 
 
\newcommand{\diver}{\operatorname{div}}
\newcommand{\dom}{G}  
\newtheorem{example}{Example} 

\title{Cross-diffusion systems with entropy structure\thanks{The author
acknowledges partial support from   
the Austrian Science Fund (FWF), grants P27352, P30000, F65, and W1245.}}

\author{Ansgar J\"ungel\thanks{Institute for Analysis and Scientific Computing, 
Vienna University of Technology, Wiedner Hauptstra\ss e 8--10, 1040 Wien, Austria 
({\tt juengel@tuwien.ac.at}).}
}

\begin{document}

\AlgLogo{1}{10}

\maketitle

\begin{abstract}
Some results on cross-diffusion systems with entropy structure are reviewed.
The focus is on local-in-time existence results for general systems with
normally elliptic diffusion operators, due to Amann, and global-in-time existence
theorems by Lepoutre, Moussa, and co-workers for cross-diffusion systems with
an additional Laplace structure. The boundedness-by-entropy method
allows for global bounded weak solutions to certain diffusion systems.
Furthermore, a partial result on the uniqueness of weak solutions is recalled,
and some open problems are presented. 
\end{abstract}



\pagestyle{myheadings}
\thispagestyle{plain}
\markboth{A. J\"UNGEL}{Cross-diffusion systems}


\section{Introduction}

Multi-species systems from physics, biology, chemistry, etc.\
can be modeled by reaction-diffusion equations. When the gradient
of the density of one species induces a flux of another species, 
cross diffusion occurs. Mathematically, this means that the diffusion
matrix involves nonvanishing off-diagonal elements.
In many applications, it turns out that the diffusion matrix is neither
symmetric nor positive definite, which considerably complicates the
mathematical analysis (see the examples in Section \ref{sec.exam} and 
\cite[Section 4.1]{Jue16}). 
In recent years, some progress has been made in
this analysis by identifying a structural condition, namely a formal
gradient-flow or entropy structure, allowing for a mathematical treatment.
In this review, we report on selected results obtained from several researchers.

The cross-diffusion equations have the form
\begin{equation}\label{1.eq}
  \pa_t u_i - \sum_{j=1}^n\diver(A_{ij}(u)\na u_j) = f_i(u)\quad\mbox{in }\Omega,\ 
	t>0,\ i=1,\ldots,n,
\end{equation}
where $u_i(x,t)$ is the density or concentration or volume fraction
of the $i$th species of a 
multicomponent mixture, $u=(u_1,\ldots,u_n)$,
$A_{ij}(u)$ are the diffusion coefficients, $f_i(u)$ is the reaction term of
the $i$th species, and $\Omega\subset\R^d$ ($d\ge 1$) is a bounded domain
with smooth boundary.
We impose no-flux and initial conditions
\begin{equation}\label{1.bic}
  \sum_{j=1}^n A_{ij}\na u_j\cdot\nu = 0\quad\mbox{on }\pa\Omega,\ t>0, \quad
	u_i(0)=u_i^0\quad\mbox{in }\Omega,\ i=1,\ldots,n,
\end{equation}
with the exterior normal unit vector $\nu$ on $\pa\Omega$,
but Dirichlet or mixed Dirichlet-Neumann boundary conditions could be considered
as well \cite{GeJu17}. Setting $A(u)=(A_{ij}(u))$ and
$f(u)=(f_1(u),\ldots,f_n(u))$, we may write \eqref{1.eq} more compactly as
$$
  \pa_t u - \diver(A(u)\na u)=f(u)\quad\mbox{in }\Omega,\ t>0.
$$

In contrast to scalar parabolic equations, generally there  
do not exist maximum principles
or a regularity theory for diffusion systems. For instance, there
exist H\"older continuous solutions to certain parabolic systems that develop 
singularities in finite time \cite{StJo95}. Here, the situation is even worse:
The diffusion matrix $A(u)$ is generally neither symmetric nor 
positive definite such that coercivity theory cannot be applied.
Our approach is to assume a structure inspired from thermodynamics:
We suppose that there exists a convex function $h:\R^n\to\R$,
called an entropy density,
such that the (possibly nonsymmetric) matrix product $h''(u)A(u)$ 
is positive semidefinite (in the sense
$z^\top h''(u)A(u)z\ge 0$ for all $z\in\R^n$). Here, $h''(u)$ denotes the
Hessian of $h$ at the point $u$.
We say that $A$ has a {\em strict} entropy structure if $h''(u)A(u)$ is positive definite
for all $u$. Then 
the entropy ${\mathcal H}[u]=\int_\Omega h(u)dx$ is a Lyapunov functional 
along solutions to \eqref{1.eq}-\eqref{1.bic} if $f(u)\cdot h'(u)\le 0$ for all
$u$:
\begin{equation}\label{1.epi}
  \frac{d{\mathcal H}}{dt} = \int_\Omega \pa_t u\cdot h'(u)dx
	= -\int_\Omega\na u:h''(u)A(u)\na udx + \int_\Omega f(u)\cdot h'(u)dx \le 0,
\end{equation}
where ``:'' denotes the Frobenius matrix product. If $h''(u)A(u)$ is positive
definite, this yields gradient estimates needed for the global existence analysis.

Introducing the entropy variables $w_i=\pa h/\pa u_i$ or $w=h'(u)$, we may
write \eqref{1.eq} equivalently as
\begin{equation}\label{1.B}
  \pa_t u(w) - \diver(B(w)\na w) = f(u(w)), \quad B(w):=A(u(w))h''(u(w))^{-1},
\end{equation}
where $u(w)=(h')^{-1}(w)$ is interpreted as a function of $w=(w_1,\ldots,w_n)$
and $h''(u)^{-1}$ is the inverse of the Hessian of $h$. By assumption,
$B(w)$ is positive semidefinite, which indicates a (nonstandard) parabolic structure. 

The entropy structure will be made more explicit for two examples in 
Section \ref{sec.exam}. In Sections \ref{sec.loc} and \ref{sec.glob}, 
the local and global in time existence of solutions, respectively, will be reviewed.
Furthermore, we comment in Section \ref{sec.more} on uniqueness results, and
we close in Section \ref{sec.open} with some open problems.


\section{Examples}\label{sec.exam}

We present two prototypic examples.

\begin{example}[Maxwell-Stefan equations]\rm\label{ex.MSE}
The dynamics of a fluid mixture of $n=3$ components with volume fractions
$u_1$, $u_2$, $u_3=1-u_1-u_2$ can be described by the
Maxwell-Stefan equations \cite{WeKr00}, defined by \eqref{1.eq} with
$$
  A(u) = \frac{1}{a(u)}\begin{pmatrix}
	d_2+(d_0-d_2)u_1 & (d_0-d_1)u_1 \\
	(d_0-d_2)u_2 & d_1+(d_0-d_1)u_2 
	\end{pmatrix},
$$
where $d_i>0$ and $a(u)=d_1d_2(1-u_1-u_2)+d_0(d_1u_1+d_2u_2)>0$.
The model can be generalized to $n\ge 3$ components; see \cite{Bot11,JuSt12}. 
For simplicity, we set $f\equiv 0$. Define the entropy density
$$
  h(u) = \sum_{i=1}^2 u_i(\log u_i-1) + (1-u_1-u_2)\big(\log(1-u_1-u_2)-1\big),
$$
where $u=(u_1,u_2)$, and the entropy ${\mathcal H}[u]=\int_\Omega h(u)dx$.
A formal computation shows that
$$
  \frac{d{\mathcal H}}{dt} 
	+ \int_\Omega\frac{1}{a(u)}\bigg(d_2\frac{|\na u_1|^2}{u_1} 
	+ d_1\frac{|\na u_2|^2}{u_2}
	+ d_0\frac{|\na(u_1+u_2)|^2}{1-u_1-u_2}\bigg)dx = 0,
$$
and in particular, $h''(u)A(u)$ is positive definite for $u_i>0$.
The entropy variables become $w_i=\pa h/\pa u_i=\log(u_i/(1-u_1-u_2))$ with
inverse $u_i(w)=e^{w_i}/(1+e^{w_1}+e^{w_2})$, which lies 
in the triangle $\dom=\{u\in\R^2:u_1,u_2>0$,
$1-u_1-u_2>0\}$. This property makes sense since $u_i$ are volume
fractions and they are expected to be bounded. 
This property can be exploited in the existence analysis to obtain 
{\em bounded} solutions without using a maximum principle (which generally cannot
be applied).
\hfill$\square$
\end{example}

\begin{example}[Population model]\rm\label{ex.SKT}
The evolution of two interacting species may be modeled by equations
\eqref{1.eq} with the diffusion matrix
$$
  A(u) = \begin{pmatrix}
	a_{10} + a_{11}u_1 + a_{12}u_2 & a_{12}u_1 \\
	a_{21}u_2 & a_{20} + a_{21}u_1 + a_{22}u_2
	\end{pmatrix},
$$
where $a_{ij}\ge 0$ \cite{SKT79}. We neglect the environmental potential
and source terms, so $f\equiv 0$. The entropy is given by
${\mathcal H}[u]=\int_\Omega h(u)dx$, where 
$h(u)=a_{21}u_1(\log u_1-1)+a_{12}u_2(\log u_2-1)$.
A formal computation shows that
\begin{equation}\label{2.pm.epi}
  \frac{d{\mathcal H}}{dt} 
	+ \int_\Omega\bigg\{\bigg(\frac{a_{10}}{u_1}+a_{21}a_{11}\bigg)|\na u_1|^2
	+ \bigg(\frac{a_{20}}{u_2}+a_{12}a_{22}\bigg)|\na u_2|^2 
	+ 4|\na\sqrt{u_1u_2}|^2\bigg\}dx
	= 0.
\end{equation}
The entropy variables are $w_1=a_{21}\log u_1$, $w_0=a_{12}\log u_2$.
Then the population densities are $u_1=e^{w_1/a_{21}}$, $u_2=e^{w_2/a_{12}}>0$.
An upper bound cannot be expected.

The model can be generalized to $n\ge 2$ species with diffusion coefficients
\begin{equation}\label{2.pm}
  A_{ij}(u) = \delta_{ij}\bigg(a_{i0} + \sum_{k=1}^n a_{ik}u_k\bigg)
	+ a_{ij}u_i, \quad i,j=1,\ldots,n.
\end{equation}
The entropy structure is more delicate than in the two-species case. Indeed,
assume that there exist numbers $\pi_i>0$ such that the equations
\begin{equation}\label{2.db}
  \pi_i a_{ij} = \pi_j a_{ji}, \quad i,j=1,\ldots,n,
\end{equation}
are satisfied. Then $h(u)=\sum_{i=1}^n\pi_i u_i(\log u_i-1)$ is an entropy density,
i.e.\ $d{\mathcal H}/dt\le 0$ \cite{CDJ17}. Equations \eqref{2.db} are
recognized as the detailed-balance condition for the Markov chain with
transition rates $a_{ij}$, and $\pi=(\pi_1,\ldots,\pi_n)$ is the corresponding
invariant measure \cite[Section 5.1]{Jue16}.
\hfill$\square$
\end{example}


\section{Local existence of classical solutions}\label{sec.loc}

A very general result on the local-in-time existence of classical solutions 
to diffusion systems was
proved by Amann (see \cite[Section 1]{Ama90} or \cite[Theorem 14.1]{Ama93}). 
A special version reads as follows.

\begin{theorem}[Amann \cite{Ama90}]
Let $\dom\subset\R^n$ be open, $A_{ij}$, $f_i\in C^\infty(\dom)$, all eigenvalues 
of $A(u)$ have positive real parts for all $u\in\dom$, and
$u^0\in V:=\{v\in W^{1,p}(\Omega;\R^n):v(\overline\Omega)\subset \dom\}$, where $p>d$.
Then there exists a unique maximal solution $u$ to \eqref{1.eq}-\eqref{1.bic}
satisfying
$u\in C^0([0,T^*);V)\cap C^\infty(\overline\Omega\times(0,T^*);\R^n)$, where
$0<T^*\le\infty$. 
\end{theorem}

An elliptic operator $u\mapsto\diver(A(u)\na u)$ with the property that
all eigenvalues of $A(u)$ have positive real parts is called {\em normally elliptic}.
We claim that any cross-diffusion system with strict entropy structure is normally
elliptic. 

\begin{lemma}[Eigenvalues of $A$]
Let $A\in\R^{n\times n}$. We assume that there exists a symmetric, positive definite
matrix $H\in\R^{n\times n}$ such that $HA$ is positive definite. Then
every eigenvalue of $A$ has a positive real part.
\end{lemma}

In the context of cross-diffusion systems, $H$ stands for the Hessian $h''(u)$.

\begin{proof}
Let $\lambda=\xi+i\eta$ with $\xi$, $\eta\in\R$ be an eigenvalue of $A$
with eigenvector $u=v+iw$, where $v$, $w\in\R^n$ with $v\neq 0$ or $w\neq 0$. 
It follows from $Au=\lambda u$
that $Av=\xi v-\eta w$, $Aw=\eta v+\xi w$. We multiply both equations by $v^\top H$,
$w^\top H$, respectively:
$$
  0 < v^\top HAv=\xi v^\top Hv - \eta v^\top Hw, \quad
	0 < w^\top HAw = \eta w^\top Hv + \xi w^\top Hw.
$$
Since $H$ is symmetric, we have $v^\top Hw=w^\top Hv$. Therefore, adding both
identities,
$$
  0 < v^\top HAv + w^\top HAw = \xi(v^\top Hv+w^\top Hw).
$$
We infer from the positive definiteness of $H$ that $\xi>0$, proving the claim.
\end{proof}


\section{Global existence of weak solutions}\label{sec.glob}

The classical solution of Amann can be continued for all time 
under some assumptions \cite[Theorem 15.3]{Ama93}.

\begin{theorem}[Amann \cite{Ama93}]
Let $u$ be the classical maximal solution to \eqref{1.eq}-\eqref{1.bic} on $[0,T^*)$.
Assume that $u|_{[0,T]}$ is bounded away from $\pa\dom$ for
each $T>0$ and that there exists $\alpha>0$ such that $\|u(t)\|_{C^{0,\alpha}}\le C(T)$ 
for all $0\le t\le T<\infty$, $t<T^*$. Then $T^*=\infty$.
\end{theorem}

Unfortunately, it is not easy to derive a uniform bound in the H\"older norm.
A possibility is to show that the gradient $\na u_i(t)$ satisfies some higher 
integrability, namely $L^p(\Omega)$ for $p>d$, since $W^{1,p}(\Omega)$ embeds
continuously into $C^{0,\alpha}(\overline\Omega)$ for $\alpha=1-p/d>0$.
Estimates in the $W^{1,p}$ norm with $p>d$ for a particular system
were derived in, e.g., \cite{HNP15,Le06}.

Another approach is to find weak solutions using the entropy method as outlined
in the introduction. The key elements of the existence proof are the definition of 
an approximate problem and a compactness argument.
We are aware of two approaches in the literature. 
In both approaches, the time derivative is
replaced by the implicit Euler discretization. This avoids issues with the
(low) time regularity. To define the change of unknowns $u(w)$, we need
bounded approximate solutions $w$. The first approach regularizes
the equations by adding a weak form of $\eps((-\Delta)^s w+w)$.
Since $H^s(\Omega)\hookrightarrow L^\infty(\Omega)$ for $s>d/2$, this yields
bounded weak solutions. The second approach formulates the implicit Euler scheme
as a fixed-point equation involving the solution operator $(M-\Delta)^{-1}$
for sufficiently large $M>0$. This allows one to exploit the regularization
property of the solution operator $(M-\Delta)^{-1}:L^p(\Omega)\to W^{2,p}(\Omega)$,
and the continuous embedding $W^{2,p}(\Omega)\hookrightarrow L^\infty(\Omega)$
for $p>d$ yields bounded solutions. 
We detail both approaches in the following subsections.

\subsection{Boundedness-by-entropy method}

This method does not only give the global existence of solutions but it
also yields $L^\infty$ bounds. 
It was first used in \cite{BDPS10} and made systematic in \cite{Jue15}.
The first key assumption is that the derivative
$h':\dom\to\R^n$ is invertible, where $\dom\subset\R^n$ is a bounded set.
Then $u(w(x,t))=(h')^{-1}(w(x,t))\in\dom$ yields lower and upper bounds for the
densities $u_i$; see Example \ref{ex.MSE}.
The second key assumption is the positive definiteness of $h''(u)A(u)$. 
Applications indicate that this property does not hold uniformly in $u$.
Therefore, we impose a weaker condition.

\renewcommand{\labelenumi}{(H\theenumi)}
\begin{enumerate}[leftmargin=9mm]
\item $h\in C^2(\dom;[0,\infty))$ is convex with invertible 
derivative $h':\dom\to\R^n$.

\item $\dom\subset(0,1)^n$ and for $z=(z_1,\ldots,z_n)^\top\in\R^n$
and $u=(u_1,\ldots,u_n)\in\dom$,
\begin{equation}\label{2.pd}
  z^\top h''(u)A(u)z \ge \kappa\sum_{i=1}^n u_i^{2m-2}z_i^2, 
	\quad\mbox{where }m\ge \frac12,
	\ \kappa>0.
\end{equation}

\item $A=(A_{ij})\in C^0(\dom;\R^{n\times n})$
and $|A_{ij}(u)|\le C_A|u_j|^a$ for all $u\in\dom$, $i,j=1,\ldots,n$, where $C_A$, $a>0$.

\item $f\in C^0(\dom;\R^n)$ and $\exists$ $C_f>0$: $\forall$ $u\in\dom$:
$f(u)\cdot h'(u)\le C_{f}(1+h(u))$.
\end{enumerate}

Hypothesis \eqref{2.pd} is satisfied with $m=\frac12$ in Examples \ref{ex.MSE} 
and \ref{ex.SKT} if $a_{10}>0$, $a_{20}>0$ and $m=1$ in Example \ref{ex.SKT} 
if $a_{11}>0$, $a_{22}>0$.
The following theorem is proved in \cite[Theorem 2]{Jue15}; 
also see \cite[Section 4.4]{Jue16}.

\begin{theorem}[Global existence \cite{Jue15}]\label{thm.ex1}
Let (H1)-(H4) hold and let $u^0\in L^1(\Omega;\R^n)$ be such that
$u^0(\Omega)\subset\overline\dom$. 
Then there exists a bounded weak solution
$u$ to \eqref{1.eq}-\eqref{1.bic} satisfying $u(\Omega,t)\subset\overline\dom$ for
all $t>0$ and $u\in L^2_{\rm loc}(0,\infty;H^1(\Omega;\R^n))$,
$\pa_t u\in L^2_{\rm loc}(0,\infty;H^1(\Omega;\R^n)')$, for all $T>0$ and
$\phi\in L^2(0,T;H^1(\Omega;\R^n))$,
$$
  \int_0^T\langle\pa_t u,\phi\rangle dt
	+ \int_0^T\int_\Omega\na\phi:A(u)\na udxdt
	= \int_0^T\int_\Omega f(u)\cdot\phi dxdt,
$$
where $\langle\cdot,\cdot\rangle$ denotes the dual pairing of $H^1(\Omega)'$,
and $u(0)=u^0$ holds in $H^1(\Omega;\R^n)'$.
\end{theorem}

The idea of the proof is to solve first for given $u^{k-1}$ the regularized problem
\begin{equation}\label{3.approx}
\begin{aligned}
  \frac{1}{\tau}\int_\Omega & \big(u(w^k)-u(w^{k-1})\big)\cdot\phi dx
	+ \int_\Omega\na\phi:B(w^k)\na w^k dx \\
  &{}+\int_\Omega\bigg(\sum_{|\alpha|=s}D^\alpha w^k\cdot D^\alpha\phi 
	+ w^k\cdot\phi\bigg)dx = \int_\Omega f(u(w^k))\cdot\phi dx
\end{aligned}
\end{equation}
for $\phi\in H^s(\Omega;\R^n)$, where $s>d/2$, $\alpha=(\alpha_1,\ldots,\alpha_n)
\in\N_0^n$ with $|\alpha|=\alpha_1+\cdots+\alpha_n=s$ is a multiindex,
$D^\alpha=\pa^s/(\pa x_1^{\alpha_1}\cdots\pa x_n^{\alpha_n})$ is a partial
derivative of order $m$, $u(w):=(h')^{-1}(w)$, and $w^k$ is an approximation
of $w(\cdot,k\tau)$ with the time step $\tau>0$. This problem is solved
by the Leray-Schauder theorem. Uniform estimates are derived from
a discrete version of the entropy-production identity \eqref{1.epi} and Hypothesis (H2). 

Let $u^{(\tau)}(x,t)=u(w^k(x))$ for $x\in\Omega$ and $t\in((k-1)\tau,k\tau]$,
$k=1,\ldots,N$, be piecewise constant functions in time.
If $t=0$, we set $u^{(\tau)}(\cdot,0)=u^0$.
We also need the time shift operator $(\sigma_\tau u^{(\tau)})(\cdot,t)=u(w^{k-1})$
for $t\in((k-1)\tau,k\tau]$. 
It follows from the boundedness and the discrete entropy-production inequality that
\cite[Section 4.4]{Jue16}
\begin{align}
  \|u^{(\tau)}\|_{L^\infty(0,T;L^1(\Omega))} &\le C, \label{3.L1} \\
  \tau^{-1}\|u^{(\tau)}-\sigma_\tau u^{(\tau)}\|_{L^2(0,T;H^s(\Omega)')}
	+ \|(u^{(\tau)})^m\|_{L^2(0,T;H^1(\Omega))} &\le C, \label{3.alc}
\end{align}
where $C>0$ is independent of $\eps$ and $\tau$.
(In fact, we have even a bound for $(u^{(\tau)})$ in $L^\infty(0,T;L^\infty(\Omega))$.)
If $m=1$, we deduce relative compactness for $(u^{(\tau)})$ in $L^2(Q_T)$
(where $Q_T=\Omega\times(0,T)$) from the discrete Aubin-Lions lemma in the
version of \cite{DrJu12}. When $m\neq 1$, we need the nonlinear version of
\cite{CDJ17,CJL14,ZaJu17}. 

\begin{lemma}[Nonlinear Aubin-Lions]\label{lem.aubin}
Let $T>0$, $m>0$, and let $(u^{(\tau)})$ be a family of nonnegative functions 
that are piecewise constant in time with uniform time step $\tau>0$. 
Assume that there exists $C>0$ such that \eqref{3.alc} holds for all $\tau>0$. 
\begin{itemize}[leftmargin=4mm]
\item Let $m>1$ and let $(u^{(\tau)})$ be bounded in $L^\infty(Q_T)$.
Then $(u^{(\tau)})$ is relatively compact in $L^p(Q_T)$ for any $p<\infty$
\cite[Lemma 9]{ZaJu17}.
\item Let $1/2 \le m\le 1$. Then $(u^{(\tau)})$ is relatively compact in 
$L^{2m}(0,T;L^{pm}(\Omega))$, where $p\ge 1/m$ and
$H^1(\Omega)\hookrightarrow L^p(\Omega)$ is continuous \cite[Theorem 3]{CJL14}.
\item Let $\max\{0,1/2-1/d\} < m < 1/2$ and let \eqref{3.L1} hold. Then
$(u^{(\tau)})$ is relatively compact in $L^1(0,T;L^{d/(d-1)}(\Omega))$
\cite[Theorem 22]{CDJ17}.
\end{itemize}
\end{lemma} 

Another version of the nonlinear Aubin-Lions lemma is shown in \cite{Mou16}. 

Theorem \ref{thm.ex1} can be directly applied to the Maxwell-Stefan equations
from Example \ref{ex.MSE} yielding the global existence of bounded weak solutions.


\subsection{Cross-diffusion system with Laplace structure}

Theorem \ref{thm.ex1} can be only applied to situations in which the densities
are bounded (volume fractions). However, the method of proof can be adapted to
cases, in which the domain $\dom$ is not bounded. The main difference
is that we cannot work in $L^\infty(\Omega)$ anymore but only in $L^p(\Omega)$
for suitable $p<\infty$. The precise value of $p$ depends on $m$ in Hypothesis
(H2), and a global existence result can be proved under certain growth conditions
on $A_{ij}(u)$ and $f_i(u)$. As an example, consider the population model
from Example \ref{ex.SKT} for $n\ge 2$ species. The following theorem was proved
in \cite{CDJ17}.

\begin{theorem}[Population model, linear $A_{ij}$ \cite{CDJ17}]\label{thm.ex2}
Let $u_i^0\ge 0$ be such that $\int_\Omega h(u^0)dx<\infty$ and let the
detailed-balance condition \eqref{2.db} and $a_{ii}>0$ hold. Then there
exists a weak solution $u=(u_1,\ldots,u_n)$ to \eqref{1.eq}-\eqref{1.bic}
with diffusion matrix \eqref{2.pm} satisfying $u_i\ge 0$ in $\Omega$, $t>0$, and
$u_i\in L^2_{\rm loc}(0,\infty;H^1(\Omega))$, $\pa_t u_i\in 
L^{q'}_{\rm loc}(0,T;W^{1,q}(\Omega)')$, where $q=2d+2$ and 
$q'=(2d+2)/(2d+1)$.
\end{theorem}

We have assumed that there is self-diffusion $a_{ii}>0$, yielding
an $L^2$ estimate for $\na u_i$, which is stronger than the $L^2$ estimate 
for $\na u_i^m$ with $m<1$. An existence result with vanishing self-diffusion
$a_{ii}=0$ was shown in \cite{ChJu06} for the two-species model. 
Here, we only have an $L^2$ bound for $\na\sqrt{u_i}$.
The lack of regularity for $\na u_i$ can be compensated by exploiting the gradient
estimate for $\na\sqrt{u_1u_2}$ in \eqref{2.pm.epi} and an $L^2\log L^2$ estimate
coming from the Lotka-Volterra reaction terms.

The detailed-balance condition can be replaced by a ``weak cross-diffusion''
assumption which is automatically satisfied if $(A_{ij})$ is symmetric; see
\cite[Formula (12)]{CDJ17}. 

Another generalization concerns {\em nonlinear} diffusion coefficients
\begin{equation}\label{3.Aij}
  A_{ij}(u) = \delta_{ij}\bigg(a_{i0} + \sum_{k=1}^n a_{ik}u_k^{s_k}\bigg)
	+ s_j a_{ij}u_iu_j^{s_j-1}, \quad i,j=1,\ldots,n,
\end{equation}
for $s_i\ge 0$.
The corresponding cross-diffusion system can be analyzed by the method of the
previous subsection. However, improved results can be obtained by exploiting
the Laplace structure, meaning that \eqref{1.eq} with \eqref{3.Aij}
writes as
\begin{equation}\label{3.Delta}
  \pa_t u_i - \Delta(u_ip_i(u)) = f_i(u), \quad\mbox{where }
	 p_i(u)=\alpha_{i0} + \sum_{j=1}^n \alpha_{ij}u_j^{s_j},
\end{equation}
and $\alpha_{ij}=a_{ij}$ for $i\neq j$ and $\alpha_{ii}=(s_i+1)a_{ii}$.
Let $a_{ii}>0$ and $s_i\le 2$. Then, by the entropy-production inequality, 
$\na u_i^{s_i/2}$ is bounded in $L^2(Q_T)$,
and the Gagliardo-Nirenberg inequality with $q=2+4/(ds_i)$ and $\theta=ds_i/(2+ds_i)$
shows that
\begin{align*}
  \|u_i^{s_i/2}\|_{L^q(Q_T)}^q &= \int_0^T\|u_i^{s_i/2}\|_{L^q(\Omega)}^q dt
	\le \int_0^T\|u_i^{s_i/2}\|_{H^1(\Omega)}^{q\theta}
	\|u_i^{s_i/2}\|_{L^{2/s_i}(\Omega)}^{q(1-\theta)}dt \\
	&\le \|u_i\|_{L^\infty(0,T;L^1(\Omega))}^{qs_i(1-\theta)/2}
	\int_0^T\|u_i^{s_i/2}\|_{H^1(\Omega)}^2 dt \le C.
\end{align*}
We deduce that $u_i$ is bounded in $L^{s_i+2/d}(Q_T)$.
Using the duality method of Pierre \cite{PiSc97}, an improved
regularity result can be derived. Indeed, set $\bar u=\sum_{i=1}^n u_i$ and
$\mu=\sum_{i=1}^n u_ip_i(u)/\bar u$. If $f_i(u)$ grows
at most linearly in $u_i$, we find that $\bar u$ solves
$\pa_t \bar u - \Delta(\mu\bar u) \le C\bar u$ 
for some constant $C>0$ depending on $f_i$.
Then (see, e.g., \cite[Lemma 1.2]{LeMo17} or the review \cite{Pie10})
\begin{equation}\label{3.dual}
  \int_0^T\int_\Omega \mu\bar u^2 dxdt \le C(T,u^0).
\end{equation}
We infer that $u_i^2p_i(u)$ is uniformly bounded in $L^1(Q_T)$, giving a bound
for $u_i$ in $L^{s_i+2}(Q_T)$. For $d>1$, 
this bound is better than the bound in $L^{s_i+2/d}(Q_T)$ derived above. 
The improved regularity is a key element in proving
the global existence of solutions \cite[Theorem 1.10]{LeMo17}
(also see the precursor versions in \cite{DLM14,DLMT15}). We define the entropy density
$h(u)=\sum_{i=1}^n h_i(u_i)$, where
$$
  h_i(u_i) = \left\{\begin{array}{ll}
	(u_i^{s_i}-s_iu_i+s_i-1)/(s_i-1) &\mbox{if }s_i\neq 1, \\
	u_i(\log u_i-1)+1 &\mbox{if }s_i=1.
	\end{array}\right.
$$

\begin{theorem}[Population model, nonlinear $A_{ij}$ \cite{LeMo17}]\label{thm.ex3}
Assume that $s_i>0$, $s_is_j\le 1$ for $i\neq j$, let the detailed-balance
condition \eqref{2.db} hold, and $f_i(u)=b_{i0}-\sum_{j=1}^n b_{ij}u_j^{\alpha_{ij}}$
for $b_{ij}\ge 0$ and $\alpha_{ij}<1$. Finally, let $u^0_i\in L^1(\Omega)\cap
H^1(\Omega)'$, $\int_\Omega h_i(u_i^0)<\infty$. Then there exists an 
integrable solution $u_i\ge 0$ to \eqref{3.Delta} and \eqref{1.bic} such that
for all smooth test functions $\phi$ satisfying $\na\phi_i\cdot\nu=0$ on $\pa\Omega$,
\begin{align*}
  -\int_0^\infty\int_\Omega & u\cdot\pa_t\phi dxdt 
	- \int_0^\infty\int_\Omega \sum_{i=1}^n u_ip_i(u)\Delta\phi_i dxdt \\
	&= \int_0^\infty\int_\Omega f(u)\cdot\phi dxdt + \int_\Omega u^0(x)\cdot\phi(x,0)dx.
\end{align*}
\end{theorem}

It is an open problem to show the same result for arbitrary $s_i>0$.

The key idea of the proof is to formulate the implicit
Euler scheme
$$
  \tau^{-1}(u_i^k-u_i^{k-1}) = \Delta F_i(u^k) + f_i(u^k), \quad\mbox{where }
	F_i(u^k)=u^k_ip_i(u^k),
$$
as the fixed-point equation
$$
  u^k = F^{-1}\Big((M-\Delta)^{-1}\big(u^{k-1}-u^k+MF(u^k)\big)\Big),
$$
where $F=(F_1,\ldots,F_n)$ and $M>0$ is a sufficiently large number.
In fact, if $M$ is large and $u_i^{k-1}>0$, 
we can show that $v:=u_i^{k-1}-u_i^k+MF_i(u_i^k) > 0$,
and by the maximum principle, $(M-\Delta)^{-1}v>0$. Then, if $F$ is a
homeomorphism on $[0,\infty)^n$, $u_i^k>0$, which yields positivity.
Moreover, elliptic regularity theory implies that for $v\in L^p(\Omega)$
with $p>d/2$, we have $(M-\Delta)^{-1}v\in W^{2,p}(\Omega)\hookrightarrow 
L^\infty(\Omega)$. This shows that $u_i^k$ is bounded in $L^\infty$
and it defines a fixed-point operator on $L^\infty(\Omega;\R^n)$.

The main assumption is that $F$ is a homeomorphism. Under this assumption,
Theorem \ref{thm.ex3} can be considerably generalized; see \cite[Theorem 1.7]{LeMo17}
for details.


\section{Uniqueness of weak solutions}\label{sec.more}

The uniqueness of weak solutions to diffusion systems is a delicate topic.
One of the first uniqueness results was shown in \cite{AlLu83}, assuming that
the elliptic operator is linear and the time derivative of $u_i$ is integrable.
The latter hypothesis was relaxed in \cite{Ott96} allowing for finite-energy
solutions but to scalar equations only. The uniqueness of solutions
was shown in \cite{PhTe17} for a cross-diffusion system with a strictly 
positive definite diffusion matrix. For cross-diffusion systems with entropy
structure (and not necessarily positive definite $A(u)$), there are much less
papers. The first result was for a special two-species population model 
\cite{JuZa16}, later extended to a volume-filling system \cite{ZaJu17},
and generalized in \cite{ChJu17} for a class of cross-diffusion systems.
In this section, we report on the result of \cite{ChJu17}.

We allow for cross-diffusion systems involving drift terms,
\begin{equation}\label{4.eq}
  \pa_t u_i = \diver\sum_{j=1}^n\big(A_{ij}(u)\na u_j+B_{ij}(u)\na\phi\big),
	\quad i=1,\ldots,n,
\end{equation}
where $\phi$ is a potential solving the Poisson equation
\begin{equation}\label{4.poisson}
  -\Delta\phi = u_0 - f(x)\quad\mbox{in }\Omega, \quad
	u_0 := \sum_{i=1}^n a_iu_i,
\end{equation}
$a_i\ge 0$ are some constants, and $f(x)$ is a given background density.
The equations are complemented by \eqref{1.bic} and $\na\phi\cdot\nu=0$
on $\pa\Omega$, $t>0$. For consistency, we need to impose the
condition $\int_\Omega\sum_{i=1}^n a_iu_i^0dx=\int_\Omega f(x)dx$.

The uniqueness proof only works for a special class of coefficients, namely
\begin{equation}\label{4.AB}
  A_{ij}(u) = p(u_0)\delta_{ij} + a_j u_iq(u_0), \quad B_{ij}(u) = r(u_0)u_i\delta_{ij},
	\quad i,j=1,\ldots,n,
\end{equation}
for some functions $p$, $q$, and $r$. The main result is as follows.

\begin{theorem}[Uniqueness of bounded weak solutions \cite{ChJu17}]\label{thm.unique}
Let $u^0\in L^\infty(\Omega)$ and $f\in L^2(\Omega)$.
Let $(u,\phi)$ be a weak solution to \eqref{4.eq}-\eqref{4.AB}, \eqref{1.bic}
such that $u_0(\Omega,t)\subset[0,L]$ for some $L>0$. 
Assume that there exists $M>0$ such that  for all $s\in[0,L]$,
\begin{align}
  & p(s)\ge 0, \quad p(s)+q(s)s\ge 0, \label{1.cond1} \\
	& r(s)s\in C^{1}([0,L]), \quad \frac{(r(s)+r'(s)s)^{2}}{p(s)+q(s)s}\le M.
	\label{1.cond2}
\end{align}
Then $(u,\phi)$ is unique in the class of solutions satisfying $\int_\Omega\phi dx=0$,
$\na\phi\in L^\infty(0,T;$ $L^\infty(\Omega))$, and
$u_i\in L^{2}(0,T;H^{1}(\Omega))$, $\pa_t u_i\in L^{2}(0,T;H^{1}(\Omega)')$
for $i=1,\ldots,n$.
In the case $r\equiv 0$, the boundedness of $u_0$ is not needed, provided that
$\sqrt{p(u_0)}\nabla u_i$, $\sqrt{|q(u_0)|}\nabla u_i\in L^2(\Omega\times(0,T))$.
\end{theorem}

The proof is based on the $H^{-1}$ method and the 
entropy method of Gajewski \cite{Gaj94}. First, we show the uniqueness of 
$u_0=\sum_{i=1}^n a_iu_i$, solving
$$
  \pa_t u_0 = \diver\big(\na Q(u_0)+R(u_0)\na\phi\big),
$$
where $Q(s)=\int_0^s(p(z)+q(z)z)dz$ and $R(s)=r(s)s$. Sine $Q$ is nondecreasing, 
the use of the $H^{-1}$ technique seems to be natural.
Given two solutions $(u,\phi)$ and $(v,\psi)$, the idea is to use the test
function $\chi$ that solves the dual problem $-\Delta\chi=u_0-v_0$ in $\Omega$,
$\na\chi\cdot\nu=0$ on $\pa\Omega$ and to show that 
$\frac{d}{dt}\|\na\chi\|_{L^2(\Omega)}^2\le C\|\na\chi\|_{L^2(\Omega)}^2$,
using the monotonicity of $Q$.
This implies that $u_0=v_0$ and $\phi=\psi$. Second, we differentiate 
(a regularized version of) the semimetric
$$
  d(u,v) = \sum_{i=1}^n\int_\Omega\bigg(h(u_i)+h(v_i)-2h\bigg(\frac{u_i+v_i}{2}\bigg)
	\bigg)dx,
$$
where $h(s)=s(\log s-1)+1$. Computing the time derivative of
$d(u(t),v(t))$, it turns out that the drift terms cancel and we end up with
$\frac{d}{dt}d(u,v)\le 0$ implying that $u=v$.

Gajewski's semimetric is related to the relative entropy or Kullback-Leibler entropy 
${\mathcal H}[u|v] = {\mathcal H}[u] - {\mathcal H}[v] 
-{\mathcal H}'[v]\cdot{\mathcal H}(u-v)$ used in statistics \cite{KuLe51}.
In fact, the proof of Theorem \ref{thm.unique} can be performed as well with the
symmetrized relative entropy $d_0(u,v)={\mathcal H}[u|v]+{\mathcal H}[v|u]$.
Both distances $d(u,v)$ and $d_0(u,v)$ behave like $|u-v|^2$ for ``small'' $|u-v|$,
but they lead to different expressions when computed explicitly.
The Kullback-Leibler entropy was also employed to derive explicit exponential
convergence rates to equilibrium \cite{CJMTU01} and to prove weak-strong
uniqueness results for (diagonal) reaction-diffusion systems \cite{Fis17}.


\section{Open problems}\label{sec.open}

We mention some open questions.

\begin{itemize}[leftmargin=4mm]
\item {\em Reaction terms:} Hypothesis (H4) excludes reaction terms which grow 
superlinearly. The global existence of solutions to cross-diffusion systems with,
for instance, quadratic reactions is an open problem. One approach could be
to consider renormalized instead of weak solutions, as done in \cite{Fis15} for
(diagonal) reaction-diffusion systems. This is currently under development
\cite{ChJu17}. Another idea is to exploit the entropy techniques devised for
reaction-diffusion systems \cite{FPT17}. 

\item {\em $n$-species population model:} It is an open problem to find global
solutions to the population model with diffusion matrix \eqref{2.pm} and $n\ge 3$
without detailed balance or ``weak cross-diffusion''. 
Numerical experiments indicate
that standard choices like the Boltzmann entropy, relative entropy, etc.\
are not Lyapunov functionals. So, the problem  to find a priori estimates is open.

\item {\em Uniqueness of solutions:} The uniqueness result presented in Theorem
\ref{thm.unique} is rather particular. One may ask whether weak-strong uniqueness
of solutions can be shown like in \cite{Fis17} for diagonal diffusion systems.
In fact, uniqueness of weak solutions is known to be delicate even for
drift-diffusion equations; see, e.g., \cite{DGJ01}.

\item {\em Regularity theory:} The duality method yields global $L^p$ regularity
results for cross-diffusion systems with Laplace structure (see \eqref{3.dual}). 
Another approach is to apply maximal $L^p$ regularity theory as done in
\cite{HMPW17} for Maxwell-Stefan systems, at least for local solutions. 
The (open) question is to what extent 
this theory can be applied to general systems with entropy structure?

\item {\em Entropies:} Given a cross-diffusion system, a major open question
is how an entropy structure can be detected. In thermodynamics, often
the entropy (more precisely: free energy) and entropy production are given
and the system of partial differential equations follows from these quantities.
Furthermore, it is an open question how large is the class of cross-diffusion
systems with entropy structure. Are there diffusion systems with
normally elliptic operator, which have no entropy structure?

\end{itemize}


\end{document}